\documentclass{amsart}

\usepackage{amsmath,amssymb,amsthm,latexsym}

\newtheorem*{theoremNN}{Theorem}
\newcommand{\btNN}{\begin{theoremNN}}
\newcommand{\etNN}{\end{theoremNN}}

\newtheorem{theorem}{Theorem}
\newcommand{\bt}{\begin{theorem}}
\newcommand{\et}{\end{theorem}}
\newtheorem{lemma}{Lemma}
\newcommand{\bl}{\begin{lemma}}
\newcommand{\el}{\end{lemma}}
\newtheorem{corollary}{Corollary}
\newcommand{\bc}{\begin{corollary}}
\newcommand{\ec}{\end{corollary}}
\newcommand{\bconj}{\begin{conjecture}}
\newcommand{\econj}{\end{conjecture}}
\newtheorem{problem}{Problem}
\newcommand{\bprob}{\begin{problem}}
\newcommand{\eprob}{\end{problem}}
\newcommand{\beq}{\begin{equation}}
\newcommand{\eeq}{\end{equation}}
\newcommand{\benum}{\begin{enumerate}}
\newcommand{\eenum}{\end{enumerate}}

\newcommand{\Z}{\ensuremath{\mathbf Z}}

\newcommand{\mcr}{\ensuremath{ \mathcal R}}

\newcommand{\bmat}{\left(\begin{matrix}}
\newcommand{\emat}{\end{matrix}\right)}
\newcommand{\bsmallmat}{\left(\begin{smallmatrix}}
\newcommand{\esmallmat}{\end{smallmatrix}\right)}

\DeclareMathOperator{\qqand}{\qquad\text{and}\qquad}

 \title[Tetrahedral differences]{Additive sumset sizes with tetrahedral differences} 
\author{Melvyn B.  Nathanson}
\address{Department of Mathematics\\Lehman College (CUNY)\\Bronx, NY 10468}
\email{melvyn.nathanson@lehman.cuny.edu}

\date{\today}

\thanks{Supported in part by  PSC-CUNY Research Award Program grant 66197-00 54.}

\begin{document}

\begin{abstract}
Experimental calculations suggest that the $h$-fold sumset 
sizes of 4-element sets of integers are concentrated at $h$ numbers 
that are differences of tetrahedral numbers.  In this paper it is proved 
that these ``popular'' sumset sizes always exist.  Explicit $h$-adically defined 
sets are constructed for each of these numbers. 
\end{abstract}

 \maketitle

\section{The sumset size problem}
The \emph{$h$-fold sum} of a set $A$ of integers, denoted $hA$, 
is the set of all sums of $h$ not necessarily distinct elements of $A$. 
If $|A| = k$, then the finite set $hA$ satisfies 
\beq          \label{h-adic:basicIneq}
h(k-1)+1 \leq |hA| \leq \binom{h+k-1}{k-1}.
\eeq
We have $|hA| = h(k-1)+1$ if and only if $A$ is an arithmetic progression of length $k$, 
and $|hA| = \binom{h+k-1}{k-1}$ if and only if $A$ is a Sidon set or a $B_h$-set, that is, a set 
such that every integer in the sumset $hA$ has a unique representation 
(up to permutation of the summands) as a sum of $h$ 
not necessarily distinct elements of $A$.  
In additive number theory, there is a huge literature on sets whose sumsets have sizes 
at the upper or lower ends of inequality~\eqref{h-adic:basicIneq}.  
For sets $A$ with $|hA|$ close to $h(k-1)+1 $, there is the theory of 
small doubling and inverse problems 
originated by Grigori Freiman (cf. Nathanson~\cite{nath96bb}). 
The survey paper of O'Bryant~\cite{obry04} reviews work on Sidon sets.  
What has been missing in additive number theory 
is the study of the full range of possible sumset sizes of 
$h$-fold sums of sets of size $k$.  

Let $\mcr_{\Z}(h,k)$ be the set of $h$-fold sumset sizes of sets of size $k$, that is, 
\[
\mcr_{\Z}(h,k) = \left\{ |hA|:A\subseteq \Z \text{ and } |A|=k \right\}.
\]
The \emph{integer interval} defined by real numbers $u$ and $v$ is the set 
\[
[u,v] = \{n \in \Z: u \leq n \leq v\}.
\]
By inequality~\eqref{h-adic:basicIneq}, 
\beq          \label{h-adic:basicSubset}
\mcr_{\Z}(h,k)  \subseteq \left[ h(k-1)+1, \binom{h+k-1}{k-1} \right].
\eeq

Not every possible sumset size is actually the size of a sumset.  
For example, relation~\eqref{h-adic:basicSubset} gives  
\[
\mcr_{\Z}(3,3)  \subseteq \left[ 7, 10 \right].
\]
We have 
\[
3\{0,1,2\} = \{0,1,2,3,4,5,6\} \qqand |3\{0,1,2\}| = 7
\]
\[
3\{0,1,3\} = \{0,1,2,3,4,5,6,7,9\} \qqand |3\{0,1,3\}| = 9
\]
\[
3\{0,1,4\} = \{0,1,2,3,4,5,6,8,9,12\} \qqand |3\{0,1,4\}| = 10 
\]
and so 
\[
 \{7,9,10\} \subseteq  \mcr_{\Z}(3,3).  
\]
However, there exists no set $A$ of integers 
with $|A| = 3$ and $|3A|=8$ (Nathanson~\cite{nath25bb}).  
Thus,  
\[
\mcr_{\Z}(3,3) = \{7,9,10\}.  
\]

This example motivates the following important problem: 
For all positive integers $h$ and $k$, compute the full range 
of sumset sizes of $h$-fold sums of sets of 
$k$ integers, that is, the set $\mcr_{\Z}(h,k)$.\footnote{There is the 
analogous problem in every additive abelian group or semigroup $G$:  Determine 
the set $\mcr_{G}(h,k)$ of the sizes of $h$-fold sums of $k$-element subsets of $G$.}
For all $h$ and $k$, we have 
\[
\mcr_{\Z}(h,1) =  \{1\} \qqand \mcr_{\Z}(1,k) =  \{k\}.
\]
Sets $A$ and $ B$ are \emph{affinely equivalent} if there exist numbers 
$\lambda \neq 0$ and $\mu$ 
such that 
\[
B = \lambda \ast A +\mu = \{\lambda a + \mu: a \in A\}. 
\]
If $A$ and $B$ are affinely equivalent, 
then $|hA| = |hB|$ for all positive integers $h$.  
Every finite set $A$ of integers is affinely equivalent to a set $B$ with $\min B = 0$ 
and $\gcd(B) = 1$.  In particular, every set of size 2 is affinely equivalent to the set $\{0,1\}$.
It follows that 
\[
\mcr_{\Z}(h,2) =  \{h+1\}. 
\]
Erd\H os and Szemer\' edi~\cite{erdo-szem83} stated that 
\[
\mcr_{\Z}(2,k) = \left[ 2k-1, \binom{k+1}{2} \right].
\]
(This is proved in~\cite{nath25bb}.)  
Thus, the unsolved problem is to determine $\mcr_{\Z}(h,k)$ 
for $h \geq 3$ and $k \geq 3$. 

A first step is to fix a positive integer $k$ and find the possible sizes 
of $h$-fold sums of sets of size $k$. 
For $k = 3$, Nathanson~\cite{nath25bb} proved that 
\[
\mcr_{\Z}(h,3) = \left\{\binom{h+2}{2} - \binom{i_0}{2} : i_0 \in [1,h] \right\}.
\]
Thus, if $|A|=3$, then $|hA|$ is a difference of triangular numbers.  
For $k =4$, the problem is still open:  
Compute
\[
\mcr_{\Z}(h,4)  \subseteq \left[ 3h+1, \binom{h+3}{3} \right].
\]
Recall that the \emph{$j$th tetrahedral number} 
$f_3^j = \binom{j+2}{3}$ is the sum of the first $j$ triangular numbers (Dickson~\cite{dick20}).  
Numerical experiments (Nathanson~\cite{nath25q} 
and O'Bryant~\cite{obry25}) suggest that, for $k = 4$, 
the ``most popular'' sumset sizes are the integers 
\[
f_3^{h+1} -  f_3^{i_0} = \binom{h+3}{3} - \binom{i_0+2}{3} 
\]
for $i_0 \in [0,h-1]$.    
These are the differences between the tetrahedral 
number  $f_3^{h+1}  = \binom{h+3}{3}$, 
which is also the size of a 4-element $B_h$-set, 
and the $h$ consecutive tetrahedral numbers 
$f_3^{0} ,  f_3^{1} ,\ldots, f_3^{h-1}$. 
It had been an open problem to decide if the  integers 
$f_3^{h+1}  - f_3^{i_0}  $ are, indeed,    
sumset sizes for all $h \geq 3$ and $i_0 \in [0,h-1]$.  
The purpose of this paper is to prove that these sumset sizes do exist 
for all positive integers $h$, that is, 
\[
\left\{\binom{h+3}{3} - \binom{i_0+2}{3} : i_0 \in [0,h-1] \right\} 
\subseteq \mcr_{\Z}(h,4)
\]
and to construct explicit $h$-adically defined sets with these sumset sizes. 

 For related work on sumset size problems in additive number theory, 
 see~\cite{fox-krav-zhan25}--\cite{schi25}.

\section{A family of $h$-adic sets} 
 
\btNN             \label{h-adic:theorem:main} 
Let $h \geq 1$.  For all $i_0 \in [0,h-1]$, let 
\[
p = 1 + (i_0-1)(h+1) 
\]
and 
\[
c = h^2 + h + 1 -p = (h+1-i_0)(h+1).  
\]
The set 
\[
A = \{0,1,h+1,c\}
\]
satisfies  $|A|=4$ and 
\[
|hA| = \binom{h+3}{3} - \binom{i_0+2}{3}.
\]
\etNN

\begin{proof} 
If $h=1$, then $i_0=1$ and $A=\{0,1,2,4\}$.  We have 
\[
|1A|=4 =  \binom{4}{3} - \binom{2}{3}.
\]

If $h \geq 2$ and $i_0 = 0$, then $p = 1 + (i_0-1)(h+1)  = -h$ and $c = h^2+h+1-p = (h+1)^2$.  
The set 
\[
A = \left\{ 0,1,h+1, (h+1)^2\right\} 
\]
is a $B_4$-set and so 
\[
|hA| = \binom{h+3}{3}  = \binom{h+3}{3}  - \binom{2}{3}. 
\]

For $h \geq 2$ and $i_0 \in [1,h-1]$,  let 
\[
B  = \{0,1,h+1\}  
\]
and 
\begin{align*}
A  = B \cup \{c\} & = \{0,1,h+1, (h+1-i_0)(h+1)\}.
\end{align*}
Note that $h+1-i_0 \geq 2$ implies $c > h+1$ and $|A| = 4$.  

We decompose the sumset $hA$ as follows:  
\[
hA  = \bigcup_{i=0}^{h} \left( (h-i)B+ic\right)  = \bigcup_{i=0}^{h} L_i 
\]
where 
\begin{align} 
L_i   & = (h-i)B+ic         \label{h-adic:Li}   \\ 
& = \bigcup_{j=0}^{h-i}\left(   (h-i-j)(h+1)  + [0,j] \right) + ic     \nonumber \\ 
& = \bigcup_{j=0}^{h-i} \left(  (h+ (h-i_0)i-j) (h+1) + [0,j] \right)  \nonumber \\
& = \bigcup_{j=0}^{h-i} M_{i,j}        \nonumber  
\end{align} 
and 
\[
M_{i,j} =  (h+ (h-i_0)i-j) (h+1) + [0,j] 
\] 
is an  integer interval  whose smallest element is a multiple  of $h+1$ 
and whose length is at most $h$.  
For $j \in [0,h-i]$, the  $h-i+1$  intervals  $M_{i,j}$ are pairwise disjoint 
and ``move to the left'' as $j$ increases from 0 to $h-i$.      
If $n \in L_i$ and $n = q(h+1) + r$ with $r \in [0,h]$, 
then $q = h+ (h-i_0)i-j$ for some $j \in [0,h-i]$ and $r \in [0,j]$, 
and so $L_i$ contains the integer interval $q(h+1) + [0, j]$. 

Because $B$ is a $B_{h}$-set with $|B|=3$, we have 
\beq            \label{h-adic:Li-size}
|L_i | = |(h-i)B+ic| = |(h-i)B| = \binom{h-i+2}{2}. 
\eeq 

For all $i \in [0,h-1]$, we have 
\begin{align*} 
\min\left(L_i\right)  & =  ic  <  (i+1)c  = \min\left( L_{i+1} \right)  
\end{align*} 
and 
\begin{align*} 
\max\left(L_i \right) & = (h-i)(h+1) + ic  \\
& < (h-i-1)(h+1) + (i+1)c \\
& = \max\left( L_{i+1} \right)  
\end{align*} 
and so the sets $L_i$ ``move to  the right'' as $i$ increases 
from 0 to $h$. 
Moreover, 
\[
\max\left( L_i \right) <  \min\left(L_{i+1} \right)  
\]
if and only if 
\[
 (h-i)(h+1)  + ic <  (i+1)c
\]
if and only if 
\[
i > h - \frac{c}{h+1} 
\]
if and only if 
\[
i \geq 1 + \left[  h - \frac{c}{h+1} \right] = i_0.  
\]
Thus, the sets $ L_i$ and $L_j$ are disjoint if $i_0 \leq i < j \leq h$ 
and, from~\eqref{h-adic:Li-size},  
\beq                    \label{h-adic:TopHalf}
\left| \bigcup_{i=i_0+1}^h L_i \right| = \sum_{i=i_0+1}^h  \left|  L_i \right| 
= \sum_{i=i_0+1}^h  \binom{h-i+2}{2}. 
\eeq
We shall prove that 
\beq                    \label{h-adic:BottomHalf}
\left| \bigcup_{i=0}^{i_0} L_i \right| 
 = \sum_{i=0}^{i_0} \binom{h-i+2}{2} - \sum_{i=0}^{i_0}  \binom{i_0+1-i}{2}. 
\eeq
Because the sets $L_i$ move to the right, we have 
\[
\left( \bigcup_{i=0}^{i_0} L_i  \right) \cap \left( \bigcup_{i=i_0+1}^h L_i \right) = \emptyset.      
\] 
Relations~\eqref{h-adic:TopHalf} and~\eqref{h-adic:BottomHalf} imply 
\begin{align*}                   \label{h-adic:PartialSum}
|hA| & = \left| \bigcup_{i=0}^h L_i \right| 
= \left|  \bigcup_{i=0}^{i_0} L_i \right| + \left|  \bigcup_{i=i_0+1}^{h} L_i \right|   \\
& = \sum_{i=0}^{i_0} \binom{h-i+2}{2} - \sum_{i=0}^{i_0}  \binom{i_0+1-i}{2} +  \sum_{i=i_0+1}^h  \binom{h-i+2}{2} \\ 
& =  \sum_{i=0}^h \binom{h-i+2}{2} - \sum_{i=0}^{i_0}  \binom{i_0+1-i}{2} \\
& =  \binom{h+3}{3} - \binom{i_0+2}{3}.  
\end{align*}
It remains to prove relation~\eqref{h-adic:BottomHalf}.  

We begin by computing $L_i \cap L_{i+t}$ for all $i \in [1,h-1]$ and $t \in [1,h-i]$.  
From relation~\eqref{h-adic:Li}, 
\[
L_i  = \bigcup_{j=0}^{h-i} \left(  (h+ (h-i_0)i-j) (h+1) + [0,j] \right) 
\]
and 
\[
L_{i+t} =  \bigcup_{j=0}^{h-i-t} \left(  (h+ (h-i_0)(i+t)-j) (h+1) + [0,j] \right)   
\]
and so $L_i \cap L_{i+t}$ is a union of intervals of the form $q(h+1)+[0,j]$ 
for integers $q$ and $j$. 
There is an integer $q$ with  $q(h+1) \in L_i \cap L_{i+t} $ if and only if 
there exist $j_0 \in [0,h-i]$ 
and $j_t \in [0,h-i-t]$ such that 
\beq           \label{h-adic:q}
q = h+ (h-i_0)i-j_0 = h+ (h-i_0)(i+t)-j_t  
\eeq
if and only if 
\begin{align*}
j_0 & = j_t - (h-i_0)t \\
&  \in [0,h-i] \cap [-(h-i_0)t, h-i - t - (h-i_0)t] \\ 
& = [0,  h-i - t -(h - i_0)t].  
\end{align*}
Conversely, if $j_0 \in  [0,  h-i - t -(h - i_0)t]$, then 
$j_t = j_0+(h-i_0)t \in [(h-i_0)t, h-i-t]$ 
and relation~\eqref{h-adic:q} is satisfied.
It follows that 
\[
q(h+1) + [0,j_0] \subseteq A_i
\]
and
\[
q(h+1) + [0,j_t] \subseteq A_{i+t}.
\]
Because $j_0 < j_t$, we have 
\[
q(h+1) + [0,j_0] \subseteq A_i \cap A_{i+t}
\]
and 
\beq           \label{h-adic:Lt-intersect}
L_i \cap L_{i+t} = \bigcup_{j_0 =0}^{h-i - t - (h - i_0)t} 
\left(  (h+ (h-i_0)i-j_0 ) (h+1) + [0,j_0 ] \right).  
\eeq
Therefore, for $i \in [0,i_0]$ and $t \in [1,h-i]$, we have 
\begin{align*}
\left|  L_i \cap L_{i+t}\right| 
& = \left|  \bigcup_{j_0=0}^{h-i - t -(h -i_0)t} 
\left(  (h+ (h-i_0)i-j_0) (h+1) + [0,j_0] \right)\right| \\ 
& =  \sum_{j_0=0}^{h-i - t -(h -i_0)t} 
\left| (h+ (h-i_0)i-j_0) (h+1) + [0,j_0] \right| \\
& = \sum_{j_0=0}^{h-i - t - (h -i_0)t} (j_0+1) \\
& = \binom{h-i - t - (h - i_0)t+2}{2}.
\end{align*}
In particular, 
\beq                  \label{h-adic:L1-intersect}
\left|  L_i \cap L_{i+1}\right| = \binom{ i_0 + 1 - i}{2}.
% i_0+1-i.
\eeq
Relation~\eqref{h-adic:Lt-intersect} also implies that, for $t \in [1,h-i]$,  
\[
L_i \setminus L_{i+t} = \bigcup_{j_0=h-i - t -(h - i_0)t +1}^{h-i}
\left(  (h+ (h-i_0)i-j_0) (h+1) + [0, j_0 ] \right) 
\]
and so 
\[
L_i \setminus  L_{i+1} \subseteq L_i \setminus  L_{i+2} \subseteq \cdots \subseteq 
L_i \setminus  L_{h}.
 \]
Therefore, 
\[
L_i \setminus \left(\bigcup_{t = 1}^{h-i} L_{i+t} \right) 
 = \bigcap_{t = 1}^{h-i} \left( L_i \setminus  L_{i+t}  \right) \\
  = L_i \setminus  L_{i+1}.
\]
The sets 
\[
L_i \setminus \left( \bigcup_{t = 1}^{h-i} L_{i+t} \right) 
\]
are pairwise disjoint for $i \in [0,h]$ and 
\[
 \bigcup_{i=0}^{i_0} L_i 
 =  \bigcup_{i=0}^{i_0} \left( L_i \setminus \bigcup_{t = 1}^{h-i} L_{i+t}  \right).
\]
Recalling~\eqref{h-adic:Li-size} and~\eqref{h-adic:L1-intersect}, we obtain  
\begin{align*}        
\left| \bigcup_{i=0}^{i_0} L_i \right| 
& = \sum_{i=0}^{i_0}  \left|   L_i \setminus \bigcup_{t = 1}^{h-i} L_{i+t} \right| 
 = \sum_{i=0}^{i_0}  \left| L_i \setminus  L_{i+1}\right|  \\ 
& = \sum_{i=0}^{i_0} \left(\left| L_i \right|  -  \left| L_i \cap L_{i+1}\right| \right)
 = \sum_{i=0}^{i_0} \left| L_i \right|  
 -  \sum_{i=0}^{i_0 }  \left| L_i \cap L_{i+1}\right| \\
& = \sum_{i=0}^{i_0} \binom{h-i+2}{2} - \sum_{i=0}^{i_0}  \binom{i_0+1-i}{2}. 
\end{align*}
This proves~\eqref{h-adic:BottomHalf} and completes the proof of the theorem. 
 \end{proof}

\section{Open problems} 
\bprob
This paper considers an important  class of $h$-adically defined 4-element sets.  
It is of interest to compute, for all $h \geq 3$ and all $p\in [0,h^2-1]$, 
the sumset sizes of the sets  
\[
A = \{0,1,h+1,h^2+h+1-p\}
\]
\eprob

\bprob
For all $h \geq 3$, compute the set of  sumset sizes of the sets 
\[
A = \{0,1,a,b\}
\]
for $2 \leq a \leq h$ and $a+1 \leq b \leq ha+1$.
\eprob

\bprob
Obtain a complete description of the sumset size set $\mcr_{\Z}(h,k)$ 
for all positive integers $h$ and $k$, 
 explain the distribution of sumset sizes for fixed $h$ and $k$, 
 and explain why some numbers cannot be sumset sizes.  
 A solution to this problem would be a fundamental theorem of additive number theory. 
\eprob

\end{document}